\documentclass[article,ij4uq]{ij4uq}              

\usepackage[hang]{footmisc}
\usepackage[ruled]{algorithm2e}
\usepackage{algorithmic}

\renewcommand{\myyear}{2018}
\renewcommand{\today}{}

\newcommand{\bA}{{\mathbb{A}}}
\newcommand{\bB}{{\mathbb{B}}}
\newcommand{\bC}{{\mathbb{C}}}
\newcommand{\bM}{{\mathbb{M}}}
\newcommand{\bE}{{\mathbb{E}}}
\newcommand{\bH}{{\mathbb{H}}}
\newcommand{\bI}{{\mathbb{I}}}

\newcommand{\bR}{{\mathbb{R}}}
\newcommand{\bN}{{\mathbb{N}}}

\newcommand{\bP}{{\mathbb{P}}}
\newcommand{\bQ}{{\mathbb{Q}}}
\newcommand{\bS}{{\mathbb{S}}}
\newcommand{\bU}{{\mathbb{U}}}
\newcommand{\bV}{{\mathbb{V}}}
\newcommand{\bW}{{\mathbb{W}}}
\newcommand{\bX}{{\mathbb{X}}}

\newcommand{\cC}{\mathcal{C}}

\newcommand{\cE}{\mathcal{E}}

\newcommand{\cL}{\mathcal{L}}

\newcommand{\cN}{\mathcal{N}}

\newcommand{\cU}{\mathcal{U}}

\newcommand{\argmax}{\operatornamewithlimits{argmax}}
\newcommand{\argmin}{\operatornamewithlimits{argmin}}

\newcommand{\bszeta}{\boldsymbol{\zeta}}

\newcommand{\bsphi}{\boldsymbol{\phi}}
\newcommand{\bspsi}{\boldsymbol{\psi}}

\newcommand{\bsvarphi}{\boldsymbol{\varphi}}

\newcommand{\bsg}{{\boldsymbol{g}}}
\newcommand{\bsp}{{\boldsymbol{p}}}
\newcommand{\bsq}{{\boldsymbol{q}}}

\newcommand{\bse}{{\boldsymbol{e}}}

\newcommand{\bsf}{\mathbf{f}}
\newcommand{\bss}{\mathbf{s}}

\newcommand{\bsv}{\mathbf{v}}
\newcommand{\bsu}{{\boldsymbol{u}}}

\newcommand{\bszero}{{\boldsymbol{0}}}

\newcommand{\beq}{\begin{equation}}
\newcommand{\eeq}{\end{equation}}
\newcommand{\ba}{\begin{array}}
\newcommand{\ea}{\end{array}}

\begin{document}

\volume{} 
\title{Hessian-based sampling for high-dimensional model reduction}
\titlehead{Hessian-based sampling for high-dimensional model reduction}
\authorhead{P. Chen,  \& O. Ghattas}
\corrauthor[1]{Peng Chen}
\author[1,2]{Omar Ghattas}
\corremail{peng@ices.utexas.edu}
\address[1]{Institute for Computational Engineering \& Sciences, The University of Texas at Austin, Austin, TX 78712.}
\address[2]{Department of Mechanical Engineering, and Department of Geological Sciences, The University of Texas at Austin, Austin, TX 78712}

\dataO{mm/dd/yyyy}
\dataF{mm/dd/yyyy}

\abstract{In this work we develop a Hessian-based sampling method for the construction of goal-oriented reduced order models with high-dimensional parameter inputs. Model reduction is known very challenging for high-dimensional parametric problems whose solutions also live in high-dimensional manifolds. However, the manifold of some quantity of interest (QoI) depending on the parametric solutions may be low-dimensional. We use the Hessian of the QoI with respect to the parameter to detect this low-dimensionality, and draw training samples by projecting the high-dimensional parameter to a low-dimensional subspace spanned by the eigenvectors of the Hessian corresponding to its dominating eigenvalues. Instead of forming the full Hessian, which is computationally intractable for a high-dimensional parameter, we employ a randomized algorithm to efficiently compute the dominating eigenpairs of the Hessian whose cost does not depend on the nominal dimension of the parameter but only on the intrinsic dimension of the QoI. We demonstrate that the Hessian-based sampling leads to much smaller errors of the reduced basis approximation for the QoI compared to a random sampling for a diffusion equation with random input obeying either uniform or Gaussian distributions. 
}


\maketitle

\section{Introduction}

Partial differential equations (PDEs) with stochastic or parametric inputs can be found in many different contexts such as uncertainty quantification, inverse problems, control and optimization, sensitivity and risk analysis. In the case that the dimension of the parameter is very high or even infinite, approximation of the parametric PDEs is computationally very challenging because of the curse of dimensionality---the computational complexity grows exponentially with respect to the dimension of the parameter. Recently, different approximation methods have been developed to deal with the high-dimensional approximation, such as Monte Carlo approximation and its variants---multilevel, quasi, high-order quasi Monte Carlo \cite{giles2015multilevel, dick2013high, dick2016higher}, sparse polynomial approximation with Galerkin projection or collocation \cite{cohen2011analytic, xiu2005high, babuvska2010stochastic}, low-rank approximation \cite{kressner2011low, nouy2015low},
Taylor approximation or perturbation analysis \cite{bonizzoni2014perturbation, AlexanderianPetraStadlerEtAl17}, and reduced basis approximation \cite{rozza2007reduced, bui2008model,  binev2011convergence, chen2016model, ChenQuarteroniRozza2017}. 

The reduced basis approximation, or more generally model reduction \cite{benner2015survey}, seek the PDE solution by a Galerkin projection in a reduced basis space that is constructed from some `snapshots'--- PDE solutions at properly selected parameter samples. The dimension of the reduced basis space is expected to be much smaller than the dimension of a high-fidelity approximation space such as the finite element space, so that solving the reduced basis system is much faster than solving the high-fidelity system. Therefore, the way to construct the reduced basis space becomes crucial for the accuracy and efficiency of the reduced basis approximation, which depends on two factors---the training samples and the construction method. For the latter, proper orthogonal decomposition (POD) or the related singular value decomposition (SVD) of the snapshot matrix, and greedy algorithms with a posteriori error estimates have been developed as the two most successful methods. For the former, the mostly used training samples are random samples drawn from the parameter space \cite{hesthaven2015certified, quarteroni2015reduced}. Quasi Monte Carlo samples \cite{quarteroni2015reduced}, structured collocation or quadrature points \cite{chen2012comparison}, sparse grid points \cite{elman2012reduced, chen2015new},  have also been used.

When it comes to problems with high-dimensional parameters, the PDE solutions may also live in high-dimensional manifolds. Therefore, a large number of reduced basis functions have to be used in order to achieve certain required accuracy of the  reduced basis approximation, which makes the reduced order model less efficient or useful. However, in many cases the goal of our computation is some QoI depending the PDE solution, e.g., the average of the solution at a certain location, which may live in a low-dimensional manifold even it depends on the high-dimensional parameter through the PDE solution. To detect this low-dimensionality structure, we use the Hessian information of the QoI with respect to the parameter, which describes its local curvature, or the extent of its variation with respect to the parameter in different directions. More specifically, one expects that the QoI varies the most along the eigenvectors corresponding to the dominating eigenvalues of the Hessian. Hence, instead of sampling in the whole parameter space, we draw samples by projecting the parameter to the subspace spanned by these eigenvectors, namely a Hessian-based sampling, which is supposed to capture the most variation of the QoI.
When the dimension of the parameter is very high, the Hessian matrix becomes very large, to form which one needs to solve a large number of PDEs that is computational intractable. To address this difficulty, we employ a randomized SVD algorithm to compute the dominating eigenpairs of the Hessian, which requires only a limited number of PDE solves. 
To demonstrate the accuracy of the Hessian-based sampling, we perform numerical experiments based on a diffusion model with parametric diffusion coefficient. We consider both a uniform distribution and a Gaussian distribution for the parameter. In the former case, the coefficient is a pieceswise random variable in each subdomain of the physical domain; for the latter, the coefficient is a log-normal random field. We construct the reduced order model by both a POD/SVD algorithm and a greedy algorithm with random training samples, as well as by the POD/SVD algorithm with Hessian-based training samples. From the comparison of the error decay of the reduced basis approximation for both the solution and the QoI, we show that the Hessian-based sampling leads to more accurate approximation for the QoI than the random sampling, not necessarily for the PDE solution. 
We mention that a Hessian-based model reduction with initial-condition inputs was developed in \cite{BashirWillcoxGhattasEtAl08}, which does not involve any parameter and the Hessian has different meaning from the second order variation of the QoI in our context. In \cite{LiebermanWillcoxGhattas2010, CuiMarzoukWillcox2016}, the parameter and state are simutaneously projected to their subspaces constructed in a greedy manner in the context of model reduction for inverse problems, and the Hessian of the likelihood function is employed in \cite{CuiMarzoukWillcox2016} to seek the parameter subspace for inverse problems.

The following of the paper is organized as follows: in Section \ref{sec:modelReduction}, we present the basic elements for model reduction, including the reduced basis approximation, offline-online decomposition, two methods for the construction of reduced basis spaces, and a short survey of sampling methods for generating the training samples. Section \ref{sec:HessianSampling} is devoted to the development of the Hessian-based sampling method, the randomized SVD algorithm for the computation of the eigenpairs of the Hessian, and the way to compute the Hessian action in certain given direction.  Numerical experiments are presented in Section \ref{sec:experiments} for the demonstration of the efficiency and the accuracy of the Hessian-based sampling method, for both a uniform distributed parameter of 256 dimensions and a Gaussian distributed parameter of 16,641 dimensions. At last, conclusions and perspectives are provided in Section \ref{sec:conclusion}. 

\section{Model reduction}
\label{sec:modelReduction}
 In this section, we briefly present the main ingredients of model reduction for a linear parametric partial differential equation (PDE) by a reduced basis method, which include a high-fidelity approximation and a reduced basis approximation for the PDE and QoI, offline-online decomposition of the reduced basis approximation, the construction algorithms (POD/SVD and greedy) of the reduced basis space, and a short survey of sampling methods for the construction. 

\subsection{Parametric partial differential equations}
Let $V$ denote a Hilbert space on $\bR$ with its dual space $V'$. Let $P \subset \bR^K$ denote a $K$-dimensional parameter space, where $K \in \bN$. We consider an abstract weak form of a linear parametric PDE: at any given parameter $\bsp = (p_1, \dots, p_K) \in P$, find $u \in V$ such that 
\beq\label{eq:weak}
a(u, v; \bsp) = f(v;\bsp) \quad \forall v \in V,
\eeq 
where $a(\cdot, \cdot; \bsp): V\times V \to \bR$ is a bilinear form and $f(\cdot;\bsp) \in V'$ is a linear functional for any given $\bsp$. By $s(u) \in \bR$ we denote a QoI that depends on the solution $u$, which is our goal of computation. 

\subsection{High-fidelity approximation}

To solve problem \eqref{eq:weak}, we introduce an approximation space $V_h \subset V$ with dimension $N_h = \text{dim}(V_h)$, e.g., a finite element space, where $h$ stands for the mesh size. In the following, we call $V_h$ a high-fidelity approximation space and $u_h$ a high-fidelity solution, which solves the high-fidelity approximation problem: at any given $\bsp \in P$, find $u_h \in V_h$ such that 
\beq\label{eq:hifi}
a(u_h, v_h; \bsp) = f(v_h; \bsp) \quad \forall v_h \in V_h.
\eeq
Let $\{\zeta_h^n\}_{n=1}^{N_h}$ denote the basis functions in $V_h$, i.e., $V_h = \text{span}\{\zeta_h^n, n = 1, \dots, N_h\}$, so that the high-fidelity solution $u_h$ can be represented as 
\beq
u_h = \sum_{n = 1}^{N_h} u_h^n \zeta_h^n,
\eeq
where $\bsu_h = (u_h^1, \dots, u_h^{N_h})^\top \in \bR^{N_h}$ is the coefficient vector. 
Then the algebraic formulation of problem \eqref{eq:hifi} can be written as: find $\bsu_h \in \bR^{N_h}$ such that  
\beq\label{eq:hifiAlg}
\bA_h(\bsp) \bsu_h = \bsf_h(\bsp),
\eeq
where the high-fidelity matrix $\bA_h(\bsp) \in \bR^{N_h \times N_h}$ and vector $\bsf_h(\bsp) \in \bR^{N_h}$ at $\bsp$ are given by 
\beq
(\bA_h(\bsp))_{mn} =  a(\zeta_h^n, \zeta_h^m;\bsp) \text{ and } (\bsf_h(\bsp))_m = f(\zeta_h^m; \bsp), \quad m, n = 1, \dots, N_h.
\eeq
As a result, the QoI $s(u)$ can be approximated by
\beq
s(u_h) = \bss_h^\top \bsu_h, \text{ where } (\bss_h)_n = s(\zeta_h^n), \quad n = 1, \dots, N_h,
\eeq
where we assume that the QoI is linear with respect to the solution for simplicity. 

\subsection{Reduced-basis approximation}
As $N_h$ is typically very big if high accuracy of the solution/QoI is required, solving the large-scale system \eqref{eq:hifiAlg} at each $\bsp \in P$ is computational expensive and only a limited number of solves can be afforded. To reduce the computational cost, we introduce a reduced basis approximation: for any given $\bsp \in P$, find $u_N \in V_N$ such that 
\beq\label{eq:rb}
a(u_N, v_N;\bsp) = f(v_N;\bsp) \quad \forall v_N \in V_N,
\eeq
where $V_N \subset V_h$ is called the reduced basis space with dimension $N$. Let $\{\zeta_N^n\}_{n=1}^N$ denote the basis functions of $V_N$, i.e., $V_N = \text{span}\{\zeta_N^n, n = 1, \dots, N\}$, then the reduced basis solution can be represented as 
\beq
u_N = \sum_{n = 1}^N u_N^n \zeta_N^n,
\eeq
with coefficient vector $\bsu_N = (u_N^1, \dots, u_N^N)^\top \in \bR^N$. Consequently, the algebraic formulation of the reduced basis approximation problem \eqref{eq:rb} can be obtained as 
\beq\label{eq:rbAlg}
\bA_N(\bsp) \bsu_N = \bsf_N(\bsp),
\eeq
 where the reduced basis matrix $\bA_N(\bsp)$ and vector $\bsf_N(\bsp)$ at $\bsp$ are given by 
 \beq
 (\bA_N(\bsp))_{mn} = a(\zeta_N^n, \zeta_N^m; \bsp) \text{ and } (\bsf_N(\bsp))_m = f(\zeta_N^m;\bsp), \quad m, n = 1, \dots, N.
 \eeq
 Moreover, the reduced basis approximation of the QoI can be evaluated as 
 \beq\label{eq:sN}
 s(u_N) = \bss^\top_N \bsu_N \text{ where } (\bss_N)_n = s(\zeta_N^n), \quad n = 1, \dots, N.
 \eeq
 \subsection{Offline-online decomposition}
Assume that the bilinear form $a$ and the linear functional $f$ allows the following affine representations with $Q_a$ and $Q_f$ terms
\beq\label{eq:affine}
a(w, v; \bsp) = \sum_{q = 1}^{Q_a} \theta_a^q(\bsp) a^q(w, v)\text{ and } f(v;\bsp) = \sum_{q = 1}^{Q_f} \theta_f^q(\bsp) f^q(v),
\eeq 
i.e., $a(\cdot, \cdot;\bsp)$ and $f(\cdot;\bsp)$ depend on the parameter $\bsp$ through the coefficients $\theta_a^q(\bsp) \in \bR$ and $\theta_f^q(\bsp) \in \bR$.
Then the reduced basis algebraic system can be written as 
\beq\label{eq:rbAlgAffine}
\left(\sum_{q = 1}^{Q_a} \theta_a^q(\bsp) \bA_N^q \right) \bsu_N = \sum_{q=1}^{Q_f} \theta_f^q(\bsp) \bsf^q_N,
\eeq
 where the reduced basis matrices $\bA_N^q \in \bR^{N\times N}$, $q = 1, \dots, Q_a$ and vectors $\bsf_N^q \in \bR^N$, $q = 1, \dots, Q_f$, are given by 
 \beq
 (\bA_N^q)_{mn} = a^q(\zeta_N^n, \zeta_N^m) \text{ and } (\bsf_N^q)_{m} = f^q(\zeta_N^m), \quad m,n = 1, \dots, N.
 \eeq
 Therefore, once the reduced basis matrices and vectors are computed and stored in the offline stage, solution of the reduced basis system \eqref{eq:rbAlgAffine} in the online stage takes $O(Q_aN^2+Q_fN)$ operations for assembling and $O(N^3)$ operations for solving, evaluation of the reduced basis approximation of the QoI takes $O(N)$ operations, which are independent of $N_h$. Thus, considerable computational reduction can be achieved by the offline-online decomposition for the solution of the parametric PDE and the evaluation of the QoI provided that $N \ll N_h$. 
 
We remark that for nonaffine or nonlinear parametric problems with possibily nonlinear QoI, an affine approximation (or so-called hyper reduction) is required to achieve an effective offline-online decomposition and computational reduction. Classical methods for such an affine approximation include empirical interpolation \cite{barrault2004empirical, maday2009general}, discrete empirical interpolation \cite{chaturantabut2010nonlinear}, weighted empirical interpolation \cite{chen2012weightedeim}, empirical operator interpolation \cite{drohmann2012reduced}, `best points' interpolation \cite{nguyen2008best}, gappy POD \cite{everson1995karhunen, bui2004aerodynamic}, GNAT \cite{carlberg2011efficient}, etc. The hyper reduction is beyond the scope of this work.
 
\subsection{Construction of the reduced basis space}\label{sec:construction}
Both the accuracy of the reduced basis approximation and the performance of the computational reduction critically depend on the reduced basis space $V_N$. 
Here we present two common algorithms for the construction of $V_N$: POD/SVD and a goal-oriented greedy algorithm. 
\subsubsection{The POD/SVD algorithm}
\label{sec:POD}
For the construction by proper orthogonal decomposition (POD), one first takes a training sample set 
\beq
\Xi_t = \{\bsp^n, n = 1, \dots, N_t\}
\eeq
with $N_t$ samples. Then the high-fidelity solution vector $\bsu_h(\bsp)$ is computed by solving the high-fidelity problem \eqref{eq:hifiAlg} at each of the training sample $\bsp \in \Xi_t$. By forming the coefficient matrix $\bU = (\bsu_h(\bsp^1), \dots, \bsu_h(\bsp^{N_t})) \in \bR^{N_h \times N_t}$, 
one then compute its singular value decomposition (SVD)
\beq
\bU = \bV \Sigma \bW^T,
\eeq
where $\bV = (\bszeta_1, \dots, \bszeta_{N_h}) \in \bR^{N_h \times N_h}$ and $\bW = (\bspsi_1, \dots, \bspsi_{N_t}) \in \bR^{N_t \times N_t}$ are orthonormal matrices, and $\Sigma = \text{diag}(\sigma_1, \dots, \sigma_{r}, 0, \dots, 0) \in \bR^{N_h \times N_t}$ is the diagonal matrix of positive singular values $\sigma_1 \geq \sigma_2 \geq \cdots \geq \sigma_r > 0$ with $r \leq \text{min}(N_h, N_t)$, denoting the rank of $\bU$. Then the reduced basis space $V_N$ is constructed with the first $N$ singular vectors of $\bV$ as the coefficient vectors for its basis functions, for $N$ such that  
\beq
N = \argmin_{n \leq r} \frac{\sum_{i = 1}^n \sigma_i^2}{\sum_{i = 1}^r \sigma_i^2} \geq 1 - \varepsilon
\eeq
with a given tolerance $\varepsilon > 0$ representing the information/energy loss. Note that the singular vectors are orthonormal in the discrete $\ell_2$-norm. To construct basis functions orthogonal with respect to the norm $X$, e.g., $L_2$-norm or energy $V$-norm, we only need to perform SVD on $\bB \bU$ where $\bX = \bB^T \bB$, e.g., a Cholesky factorization of $\bX\in \bR^{N_h \times N_h}$ where $\bX_{mn} = (\zeta_h^n, \zeta_h^m)_X$, and construct $V_N$ with $\bB^{-1} \bszeta_n$ as the coefficient vector of its $n$-th basis function $\zeta_N^n$, so that $(\zeta_N^m, \zeta_N^n)_X = \bszeta_n^T \bB^{-T} \bX \bB^{-1} \bszeta_m = \delta_{mn}$.

\subsubsection{The Greedy algorithm}
Different from the POD/SVD algorithm, the greedy algorithm seeks to construct $V_N$ iteratively in the parameter training set $\Xi_t$. At the initial step,  one often picks the first sample $\bsp^1$ from $\Xi_t$, solve the high-fidelity problem \eqref{eq:hifiAlg} at $\bsp^1$, and construct $V_1 = \text{span}\{u_h(\bsp^1)\}$. Then, for $N = 1, 2, \dots$, one chooses the next sample as 
\beq
\bsp^{N+1} = \argmax_{\bsp \in \Xi_t} \Delta_N(\bsp), 
\eeq
solve the high-fidelity problem at $\bsp^{N+1}$, and enrich $V_{N+1} = V_N \oplus \text{span}\{u_h(\bsp^{N+1})\}$, which is often orthogonalized by Gram--Schmidt process. Here $\Delta_N(\bsp)$ is an a-posteriori error indicator of the solution error $||u_h(\bsp) - u_N(\bsp)||_V$ or the goal-oriented error $|s(u_h(\bsp)) - s(u_N(\bsp))|$. 
As our goal is the computation of the QoI, we consider a goal-oriented (dual-weighted residual) error indicator for the latter, which is defined as
\beq\label{eq:errorindicator}
\Delta_N(\bsp) := |r(\psi_N;\bsp)| = |f(\psi_N;\bsp) - a(u_N, \psi_N;\bsp)|,
\eeq
where $\psi_N$ is the solution of the dual problem: given $\bsp \in P$, find $\psi_N \in W_N$ such that 
\beq\label{eq:rbdual}
a(w_N, \psi_N;\bsp ) = s(w_N) \quad \forall w_N \in W_N,
\eeq
where the reduced basis space $W_N$ can be constructed as $W_N = \text{span}\{\psi_h(\bsp^n), n = 1, \dots, N\}$ with $\psi_h(\bsp)$ the high-fidelity solution of the dual problem \eqref{eq:rbdual} in $V_h$. 
Under the affine assumption \eqref{eq:affine}, we can evaluate the weighted residual by
\beq\label{eq:errorindicatoraffine}
r(\psi_N;\bsp) = \sum_{q = 1}^{Q_f} \theta_f^q(\bsp) (\bar{\bsf}_N^q)^T \bspsi_N - \sum_{q=1}^{Q_a} \theta_a^q(\bsp) \bsu_N^T \bar{\bA}_N^q \bspsi_N,
\eeq
with $O(Q_fN +Q_aN^2)$ operations, independent of $N_h$, where 
\beq
(\bar{\bsf}_N^q)_m = f_q(\eta_N^m) \text{ and } (\bar{\bA}_N^q)_{mn} = a_q(\zeta_N^n, \eta_N^m), \quad m, n = 1, \dots, N,
\eeq
are computed and stored for once. Here by $\{\eta_N^n\}_{n=1}^N$ we denote the basis functions of $W_N$, which are obtained by Gram--Schmidt orthogonalization from $\{\psi_h(\bsp^n)\}_{n=1}^N$. 

\subsection{A short survey of sampling methods}
Both the POD/SVD construction and the greedy construction algorithms require a training sample set $\Xi_t$, which plays a crucial role in the approximation property of the reduced basis space $V_N$, especially in the case of high-dimensional parameter. On the one hand, $\Xi_t$ should be rich enough such that the main information of the manifold of the solution or the QoI can be captured by the snapshots in the training set. On the other hand, the size of $\Xi_t$ should not be redundantly large as one has to solve the expensive high-fidelity problem at each of the training sample by the POD/SVD algorithm, or compute the error indicator $\Delta_N(\bsp)$ at each of the training sample for each $N = 1, 2, \dots,$ by the greedy algorithm. 

One of the most widely used sampling method is random sampling from the probability distribution of the parameter \cite{rozza2007reduced, hesthaven2015certified}. It is rather straightforward and does not take the property of the computational QoI into account. A variant is the quasi-random sampling using low-discrepancy sequences \cite{quarteroni2015reduced, niederreiter1992random}, such as Halton or Sobol sequence, which tends to provide more equidistributed samples in the parameter space. For different probability distributions of the parameter, weighted reduced basis/POD methods \cite{chen2012weighted, spannring2017weighted, venturi2018weighted} were developed by sampling from the probability distribution with a weighted a-posteriori error estimator for the construction of the reduced basis space. 
Structured sampling methods using quadrature/collocation points such as Chebyshev points and Gauss Legendre/Hermite points have also been investigated \cite{chen2012comparison} in comparison with the random sampling methods.
In \cite{eftang2010hp}, an ``hp" adaptive sampling method was proposed, where the parameter domain is decomposed into smaller subdomains and in each subdomain a random sampling is used. An adaptive greedy sampling algorithm was proposed in \cite{hesthaven2014efficient} by adaptively cleaning and enriching the training sample set with random samples. A goal-oriented sampling method was developed in \cite{chen2013accurate} in the context of failure probability computation, where the samples are adapted to the critical limit state surface. In high-dimensional parameter space, 
a greedy sampling method in combination with the isotropic sparse grid and dimension-adaptive sparse grid has been developed in \cite{elman2012reduced} and  \cite{chen2015new}. Sampling from a subspace of the parameter space constructed using Karhunen--Lo\`eve expansion or gradient information were investigated in \cite{haasdonk2013reduced, carlberg2011low, tezzele2017combined}. 

\section{Hessian-based sampling}
\label{sec:HessianSampling}
In this section, we develop a new sampling method particularly suited for high-dimensional parametric problems based on the Hessian of the QoI with respect to the parameter. The rationale is that even the intrinsic dimension of the solution manifold is high, that of the QoI manifold could still be low, which can be captured by the low rank structure or fast spectral decay of the Hessian of the QoI with respect to the parameter. In fact, the low rank or fast spectral
decay property of the Hessian has been proven for some
specific problems and observed numerically for many others
\cite{BashirWillcoxGhattasEtAl08, FlathWilcoxAkcelikEtAl11,
  Bui-ThanhGhattas12a, Bui-ThanhGhattas13a, Bui-ThanhGhattas12,
  Bui-ThanhBursteddeGhattasEtAl12_gbfinalist,
  Bui-ThanhGhattasMartinEtAl13, ChenVillaGhattas2017,
  AlexanderianPetraStadlerEtAl16, AlexanderianPetraStadlerEtAl17,
  AlexanderianPetraStadlerEtAl14, CrestelAlexanderianStadlerEtAl17,
  PetraMartinStadlerEtAl14, IsaacPetraStadlerEtAl15,
  MartinWilcoxBursteddeEtAl12, Bui-ThanhGhattas15}.
We can therefore draw samples by projecting the high-dimensional parameter into a low-dimensional subspace spanned by the eigenvectors corresponding to the largest (absolute) eigenvalues. 

\subsection{Hessian}
\begin{figure}[!htb]
\begin{center}
\includegraphics[scale=0.5]{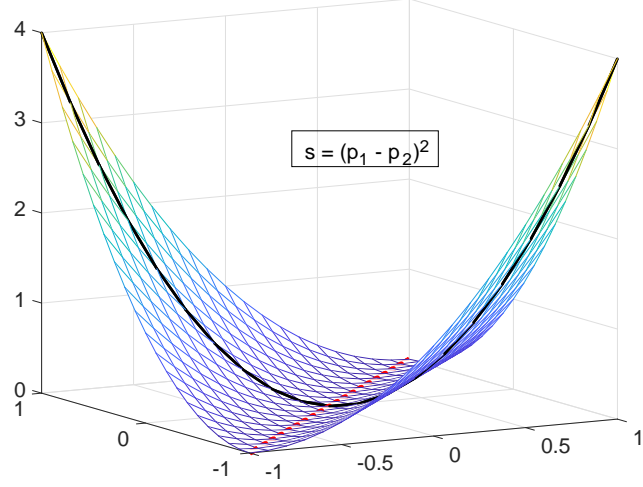}
\end{center}
\caption{The function $s(\bsp) = (p_1- p_2)^2$ with Hessian $\bH = [2, -2; -2, 2]$. Black solid line corresponding to eigenpair $\lambda_1 = 4$ and $\bsvarphi_1 = (\sqrt{2}/2,-\sqrt{2}/2)$; red dash line corresponding to eigenpair $\lambda_2 = 0$ and $\bsvarphi_2 = (\sqrt{2}/2,\sqrt{2}/2)$.}\label{fig:hessian}
\end{figure}

Our computational goal is the QoI $s(u(\bsp))$, which depends the parameter $\bsp \in P \subset \bR^{K\times K}$ through the PDE solution $u(\bsp)$. In the following, we simply denote it as $s(\bsp)$. Hessian is the square matrix $\bH \in \bR^{K\times K}$ of the second-order partial derivatives of $s$ with respect to $\bsp$, i.e., 
\beq
\bH_{kl} = \frac{\partial^2 s}{\partial p_k \partial p_l}, \quad k, l \in 1, \dots, K.
\eeq
It describes the local curvature of $s$ at $\bsp$ in the parameter space $P$, and has been widely used in large-scale optimization \cite{nocedal2006numerical, biegler2003large, ChenVillaGhattas18}, Bayesian inversion \cite{MartinWilcoxBursteddeEtAl12, biegler2011large, ChenVillaGhattas2017}, and data assimilation \cite{BashirWillcoxGhattasEtAl08, law2015data}. The eigenvectors corresponding to the dominating eigenvalues of the Hessian are the directions along which the QoI changes the most in the parameter space, which is illustrated by a simple example in Fig. \ref{fig:hessian}. We can see that the QoI varies only along the first eigenvector and does not change along the second.
Thus, sampling in the subspace spanned by the eigenvectors corresponding to the dominating eigenvalues will presumably capture the most variation of the QoI.

\vspace*{0.2cm}

\subsection{Hessian-based sampling}
\label{sec:Hessian-based sampling}

Without loss of generality, suppose the parameter $\bsp$ obeys probability distribution $\mu$ and has mean $\bar{\bsp}$ and covariance $\bC$. For instance, $\mu = \cU([-\sqrt{3}, \sqrt{3}]^K)$, a uniform distribution in the box $[-\sqrt{3}, \sqrt{3}]^K$ with mean $\bar{\bsp} = \bszero$ and covariance $\bC = \bI$, or $\mu = \cN(\bar{\bsp}, \bC)$, a Gaussian distribution with mean $\bar{\bsp}$ and covariance $\bC$. Note that a parameter $\bsp \in P$ in the compact support $P \in \bR^K$without any prescribed probability distribution can be regarded as a random variable uniformly distributed in the parameter space $P$. To proceed with the Hessian-based construction of the subspace for the parameter projection, we first consider the quadratic approximation of $s$ given by 
\beq
s(\bsp) \approx s_{\text{quad}}(\bsp) = s(\bar{\bsp}) + \bsg_{\bar{\bsp}}^T (\bsp - \bar{\bsp}) + \frac{1}{2} (\bsp-\bar{\bsp})^T \bH_{\bar{\bsp}}(\bsp-\bar{\bsp}) ,
\eeq
where $\bsg_{\bar{\bsp}}$ and $\bH_{\bar{\bsp}}$ represent the gradient and the Hessian of $s$ at $\bar{\bsp}$.
The expectation of $s$ can thus be approximated by $\bE[s] \approx \bE[s_{\text{quad}}]$, which has the analytic expression (see the proof in {Appendix A})
\beq\label{eq:meanQuad}
\bE[s_{\text{quad}}] = s(\bar{\bsp}) + \frac{1}{2} \text{tr}(\tilde{\bH}_{\bar{\bsp}}),
\eeq
where the second term is the trace of the covariance preconditioned Hessian $\tilde{\bH}_{\bar{\bsp}} = \bC \bH_{\bar{\bsp}}$ at the mean $\bar{\bsp}$. It is equivalent to the sum of all its eigenvalues, i.e.,
\beq
\text{tr}(\tilde{\bH}_{\bar{\bsp}}) = \sum_{k =1}^K\lambda_k(\tilde{\bH}_{\bar{\bsp}}). 
\eeq
In many problems, e.g.,
\cite{BashirWillcoxGhattasEtAl08, FlathWilcoxAkcelikEtAl11,
  Bui-ThanhGhattas12a, Bui-ThanhGhattas13a, Bui-ThanhGhattas12,
  Bui-ThanhBursteddeGhattasEtAl12_gbfinalist,
  Bui-ThanhGhattasMartinEtAl13, ChenVillaGhattas2017,
  AlexanderianPetraStadlerEtAl16, AlexanderianPetraStadlerEtAl17,
  AlexanderianPetraStadlerEtAl14, CrestelAlexanderianStadlerEtAl17,
  PetraMartinStadlerEtAl14, IsaacPetraStadlerEtAl15,
  MartinWilcoxBursteddeEtAl12, Bui-ThanhGhattas15}, it can be proven or numerically demonstrated that the (absolute) eigenvalues are dominated by a only a few of them $L \ll K$. Moreover, $L$ typically does not change even $K$ becomes bigger, e.g., as the mesh is refined for a random field parameter.
Therefore, the variation of $s_{\text{quad}}$ can be captured by the dominating eigenvalues, which implies that the parameter in the subspace spanned by the corresponding eigenvectors contribute to the most variation of the QoI in the parameter space. 
%

To compute the dominating eigenvalues $(\lambda_k)_{k=1}^L$ of $\tilde{\bH}_{\bar{\bsp}}$ for some $L \leq K$, which are the same as the dominating generalized eigenvalues of $(\bH_{\bar{\bsp}}, \bC^{-1})$, we solve the generalized eigenvalue problem
\beq\label{eq:gEVPmean}
{\bH}_{\bar{\bsp}} \bsvarphi_k = \lambda_k \bC^{-1}\bsvarphi_k, \text{ such that } \bsvarphi_k^T \bC^{-1} \bsvarphi_{k'} = \delta_{kk'}, \quad k, k' = 1, \dots, L.
\eeq
We remark that $\bC^{-1}$ is used in the computation as it is often readily available, e.g., when the covariance is given by the discretization of an inverse of a fractional elliptic operator as shown later in Section \ref{sec:Gaussian}.
The parameter dimension reduction is then obtained by projecting the parameter $\bsp - \bar{\bsp}$ to the $L$-dimensional subspace spanned by the eigenvectors $ \Phi_L = \text{span}\{\bsvarphi_l, l = 1, \dots, L\}$, with properly chosen $L \leq K$, i.e., 
\beq
P_L (\bsp - \bar{\bsp}) = \sum_{l =1}^L \bsvarphi_l \bsvarphi_l^T \bC^{-1}(\bsp - \bar{\bsp}).
\eeq  
Then projected parameter (sample from the subspace $\Phi_L$), denoted as $\bsp_L$, is given by 
\beq\label{eq:bspL}
\bsp_L = \bar{\bsp} + P_L(\bsp - \bar{\bsp}). 
\eeq
In the case of Gaussian distribution  $\mu = \cN(\bar{\bsp}, \bC)$, the parameter $\bsp$ can be expressed by the Karhunen--Lo\`eve expansion as 
\beq
\bsp = \bar{\bsp} + \sum_{k =1}^K\sqrt{\rho_k} \bsphi_k   \xi_k, \quad \text{i.i.d. } \xi_k \sim  \cN(0,1),
\eeq
where $(\rho_k, \bsphi_k)_{k =1}^K$ are the eigenpairs of the covariance $\bC$. Then the projection \eqref{eq:bspL} becomes 
\beq
\bsp_L =  \bar{\bsp} + \sum_{l =1}^L\bsvarphi_l \sum_{k =1}^K \sqrt{\rho_k} \bsvarphi_l^T \bC^{-1} \bsphi_k \xi_k. 
\eeq
Since a linear combination of Gaussian random variations is still a Gaussian random variable denoted as $\sqrt{\beta_l} \omega_l$, $\omega_l \in \cN(0, 1)$, with the variance given by 
\beq
\beta_l = \sum_{k =1}^K \rho_k (\bsvarphi_l^T \bC^{-1} \bsphi_k)^2 = \bsvarphi_l^T \left(\sum_{k=1}^K \rho_k^{-1} \bsphi_k \bsphi_k^T\right) \bsvarphi_l=\bsvarphi_l^T \bC^{-1} \bsvarphi_l = 1,
\eeq
where we used \eqref{eq:gEVPmean} in the last equality, so that we can sample $\bsp_L$ simply as 
\beq\label{eq:GaussSampling}
\bsp_L= \bar{\bsp}  + \sum_{l =1}^L\bsvarphi_l \omega_l.
\eeq
Note that for high-dimension parameters with $K \gg L$, $\omega_l$ can be taken as i.i.d. random variables. 

We remark that the Hessian $\bH_{\bar{\bsp}}$ is local, evaluated at the mean $\bar{\bsp}$, which may fail to characterize the variation of the QoI globally in the parameter space. To deal with this issue,
we propose two schemes for the computation of a global Hessian---namely, an averaged Hessian and a combined Hessian---to account for the variation of the QoI globally in the parameter space, as presented in \ref{app:globalHessian}. Extension of the Hessian-based sampling for multiple quantities of interest, or a vector-valued output is presented in \ref{sec:mulQoI}.

\subsection{Randomized SVD for generalized eigenvalue problems} 
\label{sec:randomizedSVD}
To solve the generalized eigenvalue problem \eqref{eq:gEVPmean},  it is prohibitive to form the full Hessian matrix when the parameter dimension is high. Instead, we apply a randomized SVD algorithm to compute the dominating generalized eigenpairs which only requires Hessian action in some random parameter directions. 
This is presented in Algorithm \ref{alg:randomizedEigenSolver}; see \cite{avron2011randomized, saibaba2016randomized} for more details on randomized SVD algorithms. 

\begin{algorithm}[!htb]
\caption{Randomized SVD for the generalized eigenvalue problem \eqref{eq:gEVPmean}}
\label{alg:randomizedEigenSolver}
\begin{algorithmic}
\STATE{\textbf{Input: } matrices $\bH_{\bsp}, \bC$, the number of eigenpairs $k$, an oversampling factor $c \leq 10$.}
\STATE{\textbf{Output: } $(\Lambda_L, \Phi_L)$ with $\Lambda_L = \text{diag}(\lambda_1, \dots, \lambda_L)$ and $\Phi_L = (\bsvarphi_1, \dots, \bsvarphi_L)$.}
\STATE{1. Draw a Gaussian random matrix $\Omega \in \bR^{n\times (L+c)}$.}
\STATE{2. Compute $Y = \bC (\bH_{\bsp} \Omega)$.}
\STATE{3. Compute $QR$-factorization $Y = QR$ such that $Q^\top \bC^{-1} Q = I_{L+c}$.}
\STATE{4. Form $T = Q^\top \bH_{\bsp} Q$ and compute eigendecomposition $T = S \Lambda S^\top$.}
\STATE{5. Extract $\Lambda_L = \Lambda(1:L, 1:L)$ and $\Phi_L = QS_L$ with $S_L = S(:,1:L)$.}
\end{algorithmic}
\end{algorithm}

We remark that the computation is dominated by the Hessian actions $\bH_{\bsp} \Omega$ and $\bH_{\bsp} Q$, which are presented in the next section. 
The advantages of Algorithm \ref{alg:randomizedEigenSolver} are: (i) the error of the eigenvalues $\lambda_j$, $j = 1, \dots, L$, are bounded by the remaining ones $\lambda_j$, $j > L$, which is small if they decay fast;
(ii) the computational cost is dominated by $2(L+c)$ Hessian actions (the application of $\bC$ on a vector is inexpensive, e.g., it only takes $O(n)$ operations by a multigrid solver for $\bC$ discretized from a differential operator); (iii) it is tractable as $L$ typically does not change when $K$ becomes bigger;
 (iv) computing the Hessian actions $\bH_\bsp\Omega$ and $\bH_\bsp Q$ can be asynchronously parallelized.

\subsection{Hessian action in a given direction}
\label{sec:HessianAction}
To compute the Hessian action in a certain given parameter direction,
we employ a Lagrange multiplier method. We take the Hessian $\bH_{\bsp}$ at any given $\bsp \in P$ as an example to compute its action in a given direction $\hat{p}\in P$, which readily applies to the local, the averaged, and the combined Hessian actions. We first form the Lagrangian as 
\beq
\cL(u, v, \bsp) = s(u) + f(v;\bsp) - a(u,v;\bsp),
\eeq 
where $v$ is the adjoint variable or the Lagrange multiplier. Then by setting the first variation of $\cL$ with respect to the adjoint and the state variables to be zero we obtain the state and adjoint variables at $\bsp$ as the solutions of the state problem \eqref{eq:weak} and the adjoint problem: find $v \in V$ such that 
\beq
a(w, v;\bsp) = s(w) \quad \forall w \in V.
\eeq
Then we can evaluate the gradient of $s$ with respect to $\bsp$ as 
\beq\label{eq:gradient}
\bsg_\bsp = \partial_\bsp \cL(u, v, \bsp) = \partial_\bsp f(v; \bsp) - \partial_\bsp a(u,v;\bsp).
\eeq
To compute the Hessian action in direction $\hat{p}$, we form another Lagrangian for the first order variation constraints as 
\beq\label{eq:LH}
\cL^H(u, v, \bsp, \hat{u}, \hat{v}, \hat{\bsp}) = a(u,\hat{v};\bsp) + a(\hat{u}, v; \bsp) + (\partial_\bsp f(v; \bsp) - \partial_\bsp a(u, v; \bsp), \hat{\bsp}),
\eeq
where $\hat{v}$ and $ \hat{u}$ are adjoint variables. They can be obtained by setting variation of \eqref{eq:LH} with respect to $u$ and $v$ as zero,
which leads to the incremental adjoint problem: find $\hat{v} \in V$ such that 
\beq\label{eq:incadj}
a(\tilde{u}, \hat{v};\bsp) =  (\partial_\bsp a(\tilde{u}, v; \bsp),  \hat{\bsp}) \quad \forall \tilde{u} \in V,
\eeq
and the incremental state problem: find $\hat{u}\in V$ such that
\beq\label{eq:incsta}
a(\hat{u}, \tilde{v};\bsp) =  (\partial_\bsp a(u, \tilde{v};\bsp) - \partial_\bsp f(\tilde{v};\bsp), \hat{\bsp}) \quad \forall \tilde{v} \in V.
\eeq
We remark that the adjoint problem and the two incremental problems are linear and have the same linear operator (or its adjoint), even when the state problem is nonlinear. Then the Hessian action in direction $\hat{\bsp}$ is given by
\beq\label{eq:HessionAction}
\bH_{\bsp} \hat{\bsp} = \partial_\bsp \cL^H \hat{\bsp} = (\partial_\bsp a(u, \hat{v};\bsp) + \partial_\bsp a(\hat{u}, v; \bsp) + \partial_{\bsp\bsp} f (v;\bsp) - \partial_{\bsp\bsp} a(u,v;
\bsp), \hat{\bsp}).
\eeq
Therefore, once the solutions of the state and adjoint problems at $\bsp$ are obtained, the Hessian action $\bH_\bsp \hat{\bsp}$ only involves solving two linear incremental problems for each $\hat{\bsp}$, which implies that the computational cost of the randomized SVD Algorithm \ref{alg:randomizedEigenSolver} is dominated by $4(L+c)$ linear PDE solves with the same linear operator (or its adjoint).

\section{Numerical experiments}
\label{sec:experiments}
To demonstrate the performance of the Hessian-based sampling algorithm for goal-oriented model reduction with high-dimensional parameter, we consider the diffusion problem 
\beq\label{eq:diffusion}
-\nabla( \kappa(\bsp) \nabla u) = g, \quad \text{ in } D,
\eeq
in a physical domain $D = [0, 1]^2$, with parametric coefficient $\kappa(\bsp)$ and suitable boundary conditions on $\partial D$. We consider the following QoI 
\beq
s(\bsp) = \frac{1}{|D_s|}\int_{D_s} u(\bsp) dx,
\eeq
where we set $D_s = [0, 0.1]^2$ with volume $|D_s| = 0.1^2$. As for the parameter $\bsp \in P \subset \bR^K$, we consider the cases of a uniform distribution with dimension $K = 16^2$ and a Gaussian distribution with dimension $K = 129^2$.

\subsection{Uniform distribution}
\label{sec:uniform}
In this example, we consider the coefficient $\kappa(\bsp)$ as a piecewise random variable with uniform distribution given by
\beq
\kappa(\bsp) = \kappa_0 + \sum_{k = 1}^K k^{-\beta} \chi_{D_k} p_k,
\eeq
where $\chi_{D_k}$ is a characteristic function taking value one in $D_k$ and zero elsewhere, $\bsp \sim \cU([-\sqrt{3}, \sqrt{3}]^K)$ with mean $\bar{\bsp} = \bszero$ and covariance $\bC = \bI$, $\beta$ is a scaling parameter.
Here $D = \cup_k D_k = [0, 1]^2$, where the $k$-th subdomain $D_k = [i*h, j*h] \times [(i+1)*h, (j+1)*h]$ with $h = 1/\sqrt{K}$, $j = \text{mod}(k, \sqrt{K})$ and $i = (k - j)/\sqrt{K}$. We take
$g = 0$, and set the Dirichlet boundary conditions $u = 1$ on $x = [0,1]\times 0$, and $u = 0$ on $x = [0,1]\times 1$, while setting zero Neumann boundary conditions elsewhere.

The weak form of problem \eqref{eq:diffusion} can be expressed as \eqref{eq:weak} with affine representation \eqref{eq:affine}, where we have $\theta_a^1(\bsp) = \kappa_0, a^1(w,v) = \int_{D} \nabla w \cdot \nabla v dx $ and 
\beq
\theta_a^{k+1}(\bsp) = k^{-\beta} \chi_{D_k} p_k, \text{ and } a^{k+1}(w, v) = \int_{D_k} \nabla w \cdot \nabla v dx, \quad k = 1, \dots, K,
\eeq
and $\theta_f^1(\bsp) = \kappa_0, f^1(v) = \int_{\Gamma_D} \nabla u_D \cdot \nabla v dx $, being $u_D$ the Dirichlet boundary condition on $\Gamma_D = [0, 1]\times \{0, 1\}$, and 
\beq
\theta_f^{k+1} (\bsp) = k^{-\beta} \chi_{D_k} p_k, \text{ and } f^{k+1}(v) = \int_{\partial D_k \cap \Gamma_D} \nabla u_D \cdot \nabla v dx, \quad k = 1, \dots, K.
\eeq
In the numerical test, we use piecewise linear finite element in a uniform mesh of size $65\times 65$ for the discretization of the problem. We consider a relatively high dimension $K = 256$. We set the parameter $\kappa_0 = \sqrt{3} + 0.01$ to guarantee that the coefficient is positive, and set $\beta = 1$ so that the solution manifold is relatively high-dimensional yet its reduced basis approximation error still shows evident decay with respect to the number of reduced basis functions. 

\begin{figure}[!htb]
\begin{center}
\includegraphics[scale=0.5]{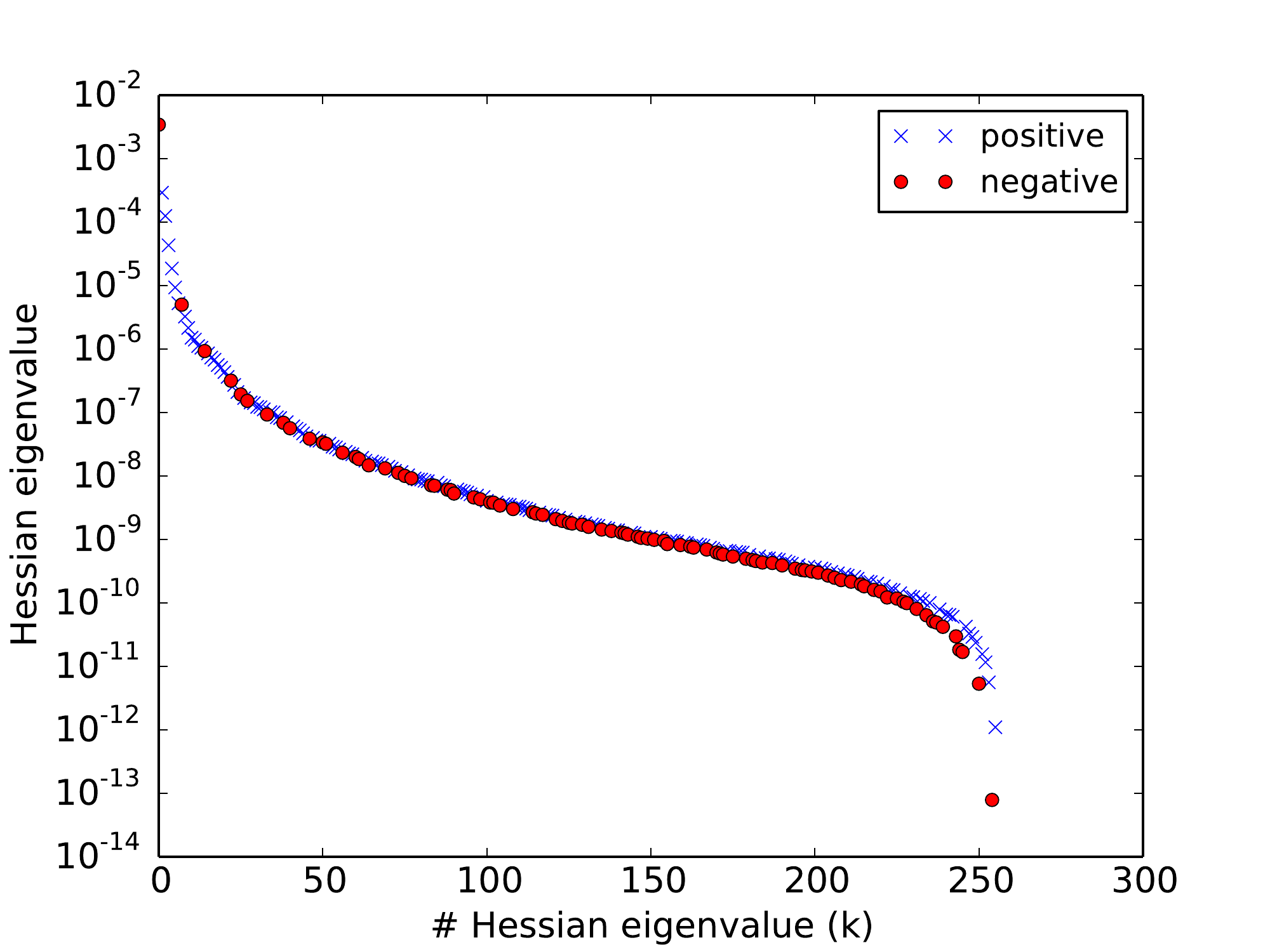}
\end{center}
\caption{The decay of the generalized eigenvalues of the Hessian $\bH_{\bar{\bsp}}$ at the mean $\bar{\bsp} = \bszero \in \bR^{256}$, as the solution of problem \eqref{eq:gEVPmean}. Positive and negative eigenvalues are shown with red dot and blue cross, respectively.}
\label{fig:dHessianMean}
\end{figure}

At first,  we compute the Hessian $\bH_{\bar{\bsp}}$ of the QoI $s$ at the mean $\bar{\bsp} = \bszero$. Note that here $\bH_{\bar{\bsp}} \in \bR^{K\times K}$, which can be formed via \eqref{eq:HessionAction} by solving the incremental adjoint and state problems \eqref{eq:incadj} and \eqref{eq:incsta} with $\hat{\bsp} = \bse_k$, whose $k$-th element is one and all the other elements are zero, $k = 1, \dots, K$.  Computing the full Hessian is used in this test of uniform distribution for $K = 256$, which becomes very expensive if $K$ is much larger as in the next test of Gaussian distribution where we use the randomized SVD Algorithm \ref{alg:randomizedEigenSolver} instead of computing the full Hessian. The eigenvalues of the Hessian are computed  
as the solution of problem \eqref{eq:gEVPmean} with $\bC = \bI$, which are shown in Fig. \ref{fig:dHessianMean}. We can observe that the eigenvalues decay very fast in the first few dimensions, with four orders of magnitude of difference in the first 20 dimensions. 

\begin{figure}[!htb]
\begin{center}
\includegraphics[scale=0.4]{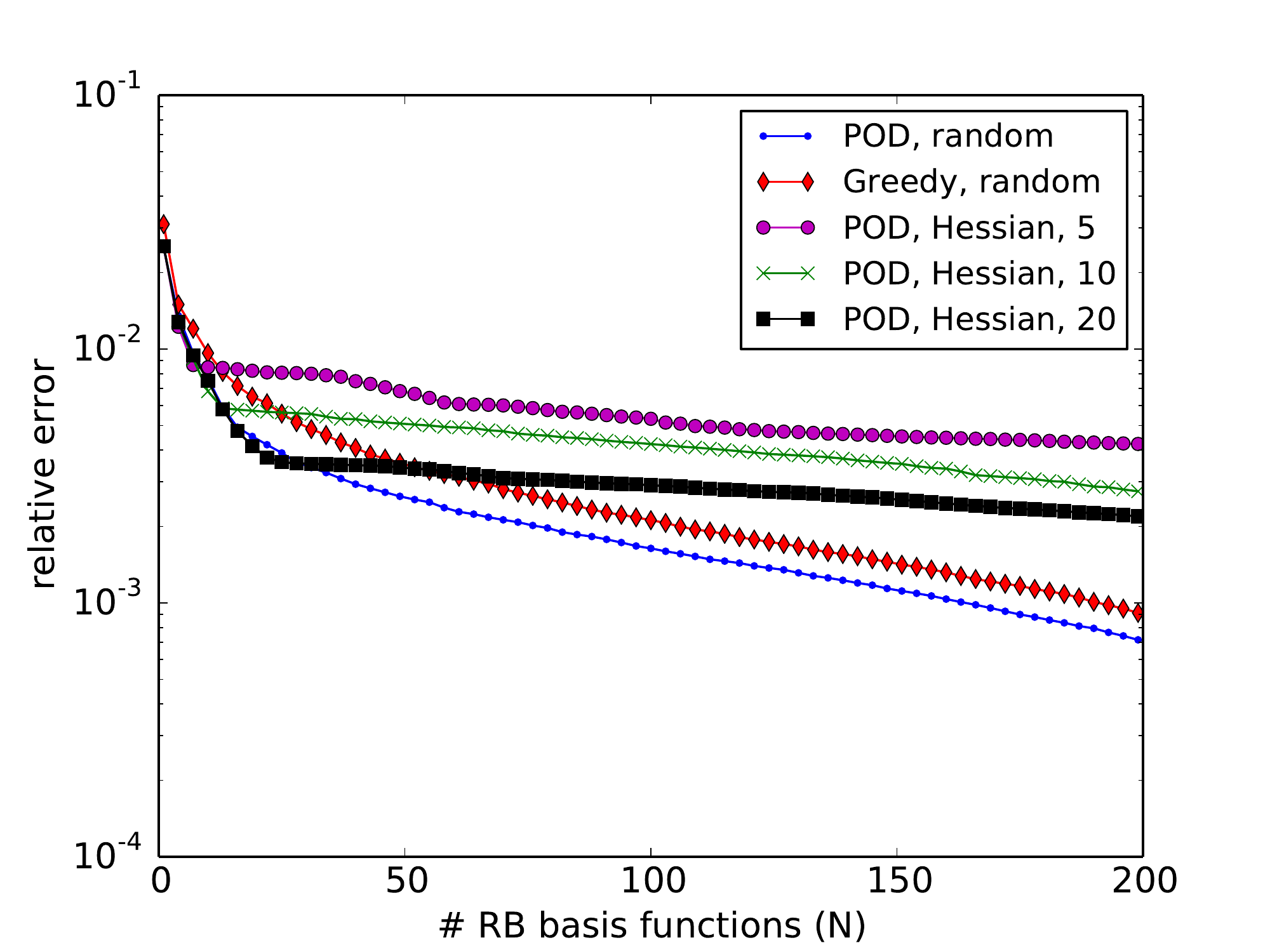}
\includegraphics[scale=0.4]{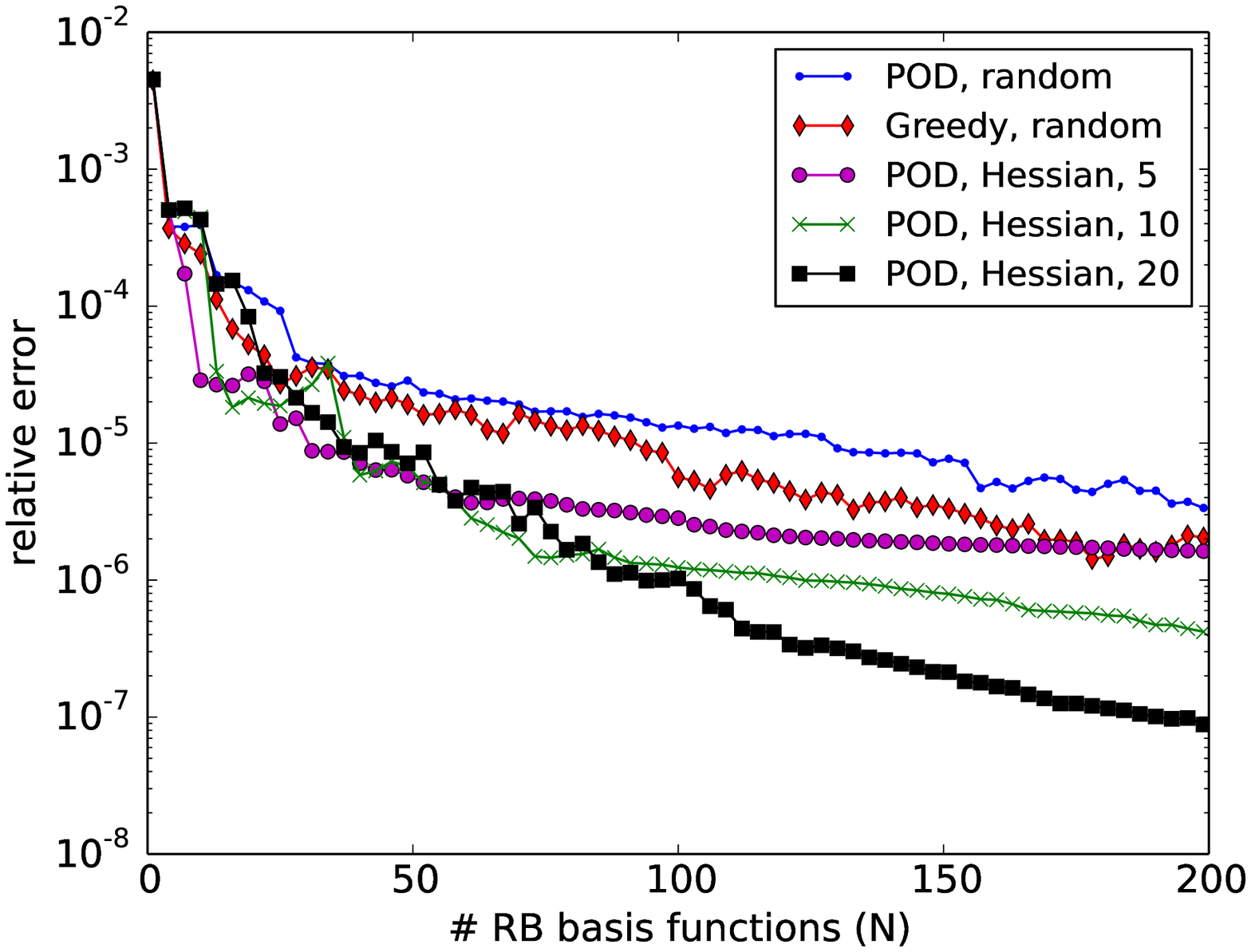}
\end{center}
\caption{Decay of the relative error $\cE_N^u$ (left) and $\cE_N^s$ (right), defined in \eqref{eq:error}, of the RB approximation of the solution and the QoI, respectively. $L = 5, 10, 20$ modes are used in the projection \eqref{eq:bspL} for the Hessian-based samples.}\label{fig:POD_randomVShessian}
\end{figure}

In the construction of the reduced basis space, we compare three different schemes: POD with random training samples, (goal-oriented) greedy with random training samples, and POD with Hessian-based training samples. We draw $N_t = 1000$ training samples for each scheme and construct the corresponding reduced basis space by the POD/greedy algorithms presented in Sec. \ref{sec:construction} with $N = 200$ basis functions. Then we compute the relative error of the reduced basis approximation for the solution and the QoI as
\beq\label{eq:error}
\cE^u_N = \frac{1}{|\Xi_{\text{test}}|}\sum_{\bsp \in \Xi_{\text{test}}}\frac{||u_h(\bsp) - u_N(\bsp)||_V}{||u_h(\bsp)||_V}; \; \cE^s_N = \frac{1}{|\Xi_{\text{test}}|}\sum_{\bsp \in \Xi_{\text{test}}}\frac{|s_h(\bsp) - s_N(\bsp)|}{|s_h(\bsp)|},
\eeq
where $\Xi_{\text{test}}$ is a test sample set with $|\Xi_{\text{test}}| = 10$ samples randomly drawn from the uniform distribution. The decay of the relative error of the RB approximation for the PDE solution is shown in the left of Fig. \ref{fig:POD_randomVShessian}, from which we can observe that the RB space obtained by POD construction with random training samples leads to the smallest approximation error, smaller than that by the greedy construction with the same random training samples. The errors get stagnated for the POD construction with the Hessian-based training samples. Moreover, the larger the number of the modes $L$ used in the projection \eqref{eq:bspL}, the smaller the errors become. This observation indicates that for the RB approximation of the solution, the Hessian-based sampling does not capture the solution manifold as well as the random sampling, which is expected since the Hessian is for the QoI, not for the PDE solution. 

As for the RB approximation of the QoI, from the right of Fig. \ref{fig:POD_randomVShessian} we can observe that the Hessian-based sampling (with $L = 20$ modes) leads to much smaller errors compared to both the POD and the greedy construction with random training samples. Moreover, the Hessian-based sampling with small number of modes ($L = 5$) yields smaller errors for small $N$ but larger errors for large $N$ than that with large number of modes ($L = 20$), which implies that the parameter in the subspace spanned by fewer eigen-modes can capture more representative reduced basis functions for small $N$, while more eigen-modes are needed if higher approximation accuracy is required. Therefore, an adaptive Hessian-based sampling by gradually increasing the number of modes could probably lead to a better construction method, which is subject to further investigation.
Note that the greedy construction yields smaller errors than the POD construction with the same random training samples, due to the use of the goal-oriented a-posteriori error indicator \eqref{eq:errorindicator} that can be efficiently computed as in \eqref{eq:errorindicatoraffine} thanks to the affine representation \eqref{eq:affine}. 

\subsection{Gaussian distribution}
\label{sec:Gaussian}
In the second example, we consider the coefficient as a log-normal random field, i.e., 
\beq
\kappa(p) = e^p,
\eeq
where $p$ is a Gaussian random field with distribution $\cN(\bar{p}, \cC)$. Here the covariance operator $\cC$ is taken as the inverse of a fractional elliptic operator
\beq
\cC = (-\delta \triangle + \gamma I)^{-\alpha}.
\eeq 
In this example, we take $\gamma = 0.5$, $\delta = 1$, $\alpha = 2$, such that $\cC$ is of trace class. By piecewise linear finite element discretization in a uniform mesh of size $129\times 129$, we obtain a $129^2$-dimensional parameter $\bsp$ as the coefficient for the random field $p$, which obeys the Gaussian distribution $\bsp \sim \cN(\bar{\bsp}, \bC)$ with covariance given such that 
\beq
\bC^{-1} = \bA \bM^{-1} \bA,
\eeq
where $\bM$ is the mass matrix and $\bA$ is the stiffness matrix given by 
\beq
\bM_{ij} = \int_D \psi_j \psi_i dx \text{ and } \bA_{ij} =  \int_D \big( \delta \nabla \psi_j \cdot \nabla \psi_i + \gamma \psi_j \psi_i \big) dx, \; i, j = 1, \dots, 129^2,
\eeq
where $\psi_i$, $i = 1, \dots, 129^2$, are the finite element basis functions. For simplicity, we take the source term $g = 1$ and use homogeneous Dirichlet boundary conditions.

\begin{figure}[!htb]
\begin{center}
\includegraphics[scale=0.5]{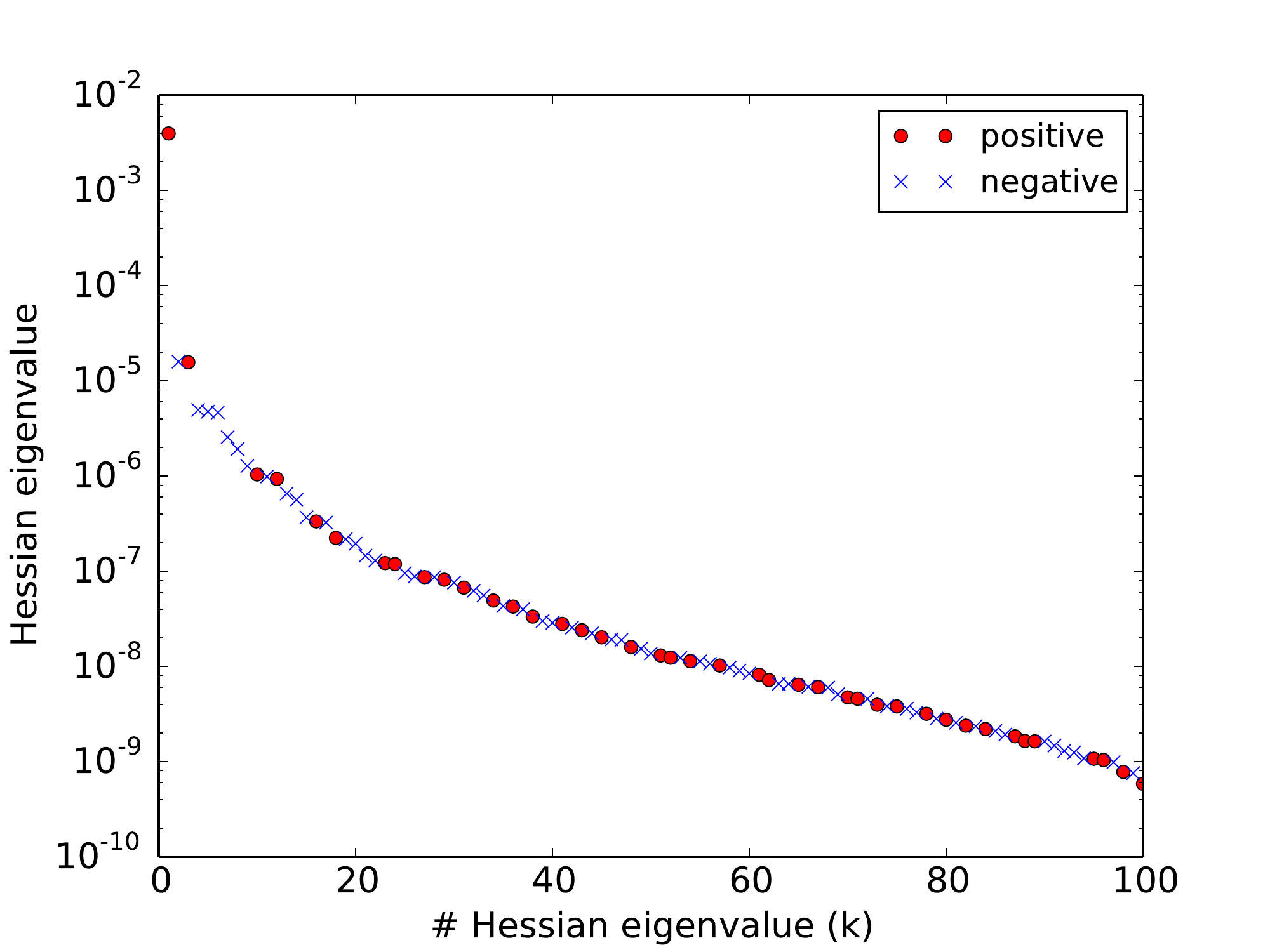}
\end{center}
\caption{The decay of the generalized eigenvalues of the Hessian $\bH_{\bar{\bsp}}$ at the mean $\bar{\bsp} = \bszero \in \bR^{16641}$, as the solution of problem \eqref{eq:gEVPmean}. Positive and negative eigenvalues are shown with red dot and blue cross, respectively.}
\label{fig:dHessianMeanPoisson}
\end{figure}

We solve the generalized eigenvalue problem \eqref{eq:gEVPmean} by the randomized SVD algorithm presented in Section \ref{sec:randomizedSVD}, where the Hessian action in a random direction is evaluated as in Section \ref{sec:HessianAction}. We compute the first 100 eigenpairs by using 110 random directions, for which the decay of the eigenvalues is shown in Fig. \ref{fig:dHessianMeanPoisson}. We can observe that the eigenvalues decay very fast in the first few dimensions, with four orders of magnitude of difference in the first 15 dimensions, and in particular the first eigenvalue is dominating, which indicates that sampling in a rather low-dimensional subspace spanned by the eigenvectors could be sufficient to capture the major variation of the QoI in the parameter space. 

\begin{figure}[!htb]
\begin{center}
\includegraphics[scale=0.4]{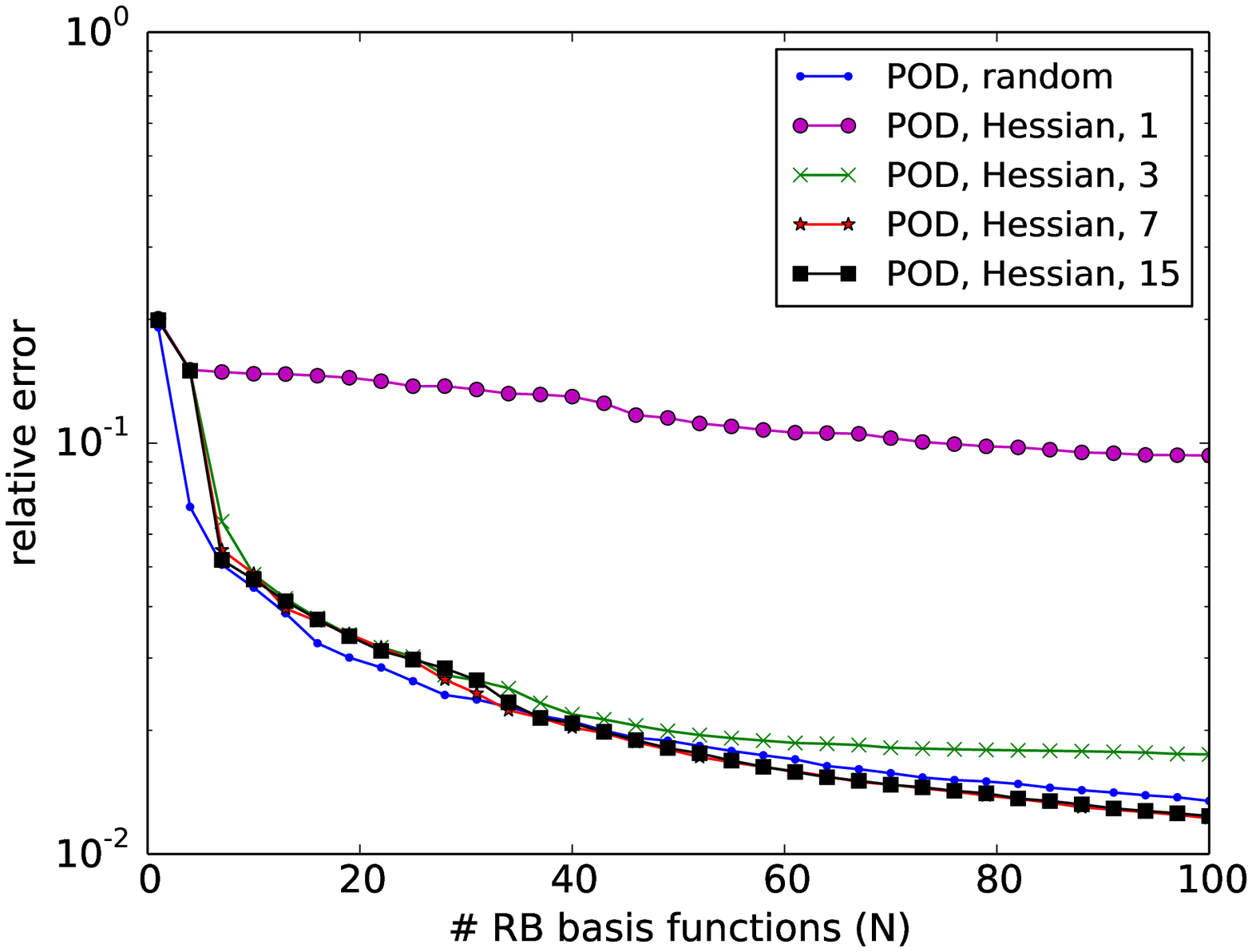}
\includegraphics[scale=0.4]{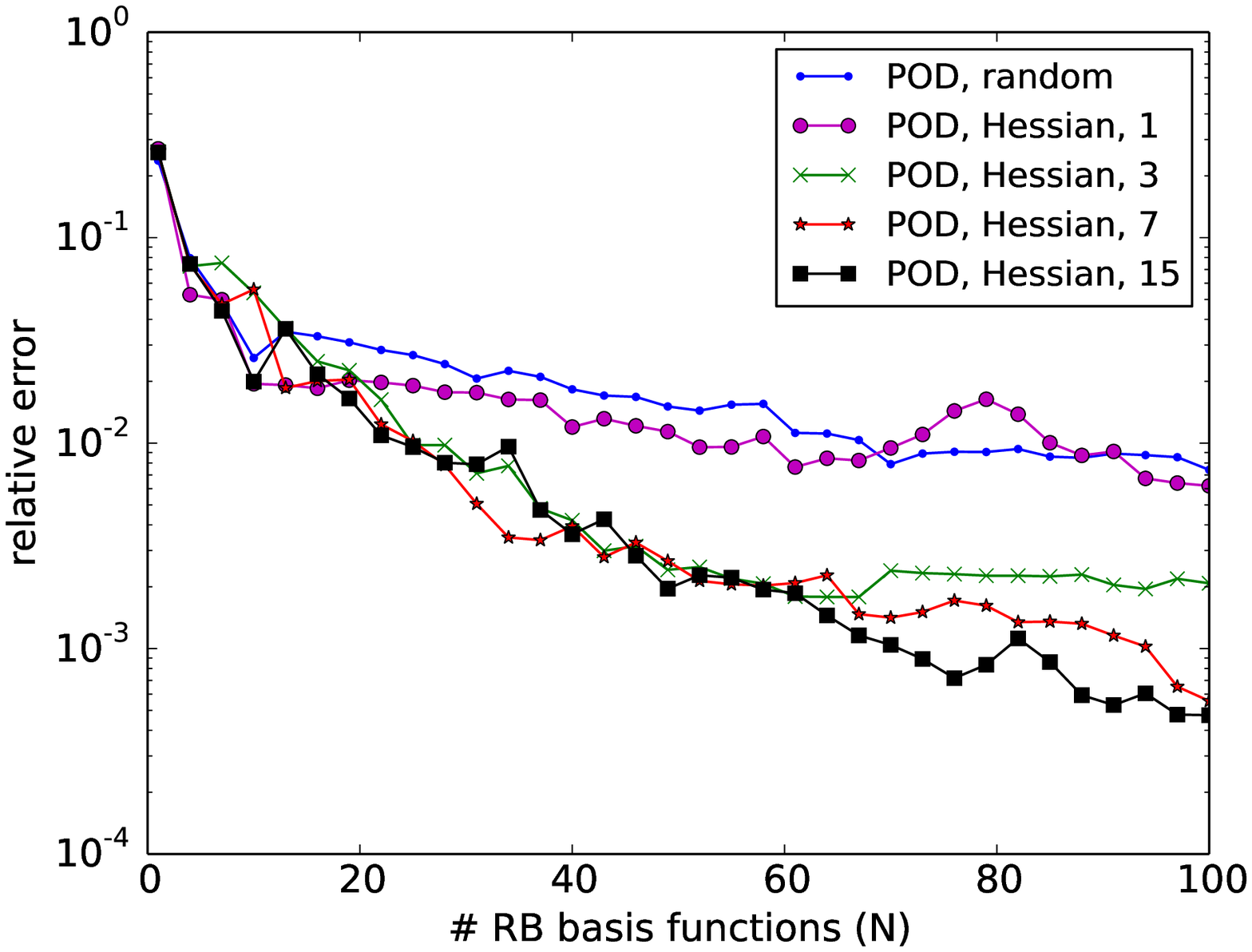}
\end{center}
\caption{Decay of the relative error $\cE_N^u$ (left) and $\cE_N^s$ (right), defined in \eqref{eq:error}, of the RB approximation of the solution and the QoI, respectively. $L = 1, 3, 7, 15$ modes are used in the projection \eqref{eq:bspL} for the Hessian-based samples.}\label{fig:POD_randomVShessian_Poisson}
\end{figure}

To demonstrate the efficacy of the Hessian-based sampling, we project the $129^2$-dimensional parameter $\bsp$ to low-dimensional subspaces as in \eqref{eq:bspL} with $L = 1, 3, 7, 15$, where the sampling is performed efficiently as in \eqref{eq:GaussSampling}. We construct the reduced basis space by the POD algorithm in Section \ref{sec:POD} with both random training samples and the Hessian-based training samples of size 1000. Note that the problem is nonaffine due to the log-normal coefficient $\kappa(\bsp) = e^\bsp$, we do not use greedy algorithm which is very expansive without the offline-online decomposition for the evaluation of the dual-weighted residual \eqref{eq:errorindicator}. Affine approximation (e.g., by empirical interpolation) of the log-normal random field is not considered here and can be found in \cite{Chen2016} for details. A Hessian-based sampling for empirical interpolation is out of scope of this paper and subject to further investigation. We compute the reduced basis approximation errors for the solution and the QoI defined in \eqref{eq:error} with 10 test samples randomly drawn from the whole parameter space.  
The decay of the errors are shown in the left of Fig. \ref{fig:POD_randomVShessian_Poisson} for the solution and in the right of Fig. \ref{fig:POD_randomVShessian_Poisson} for the QoI. 
From the former figure we can see that with only 1 mode for the Hessian-based sampling, the RB error for the solution remains large, while with 3 modes, the RB error becomes much smaller, and with 7 and 15 modes, the RB errors are comparable to and become even smaller than that obtained by random samples at large number of RB basis functions. On the other hand, with 1 mode for the Hessian-based sampling, the RB error for the QoI is already close to that obtained by random samples as seen from the right of Fig. \ref{fig:POD_randomVShessian_Poisson}. Moreover, with 3, 7, and 15 modes, the RB errors become much smaller than that obtained by random samples, which demonstrate the efficiency of the Hessian-based sampling in capturing the QoI variation in high-dimensional parameter space.

\section{Conclusion}
\label{sec:conclusion}
We developed an efficient Hessian-based sampling method to construct goal-oriented reduced order models for high-dimensional parametric problems. Even the dimension of the solution manifold is high due to the high-dimensionality of the parameter space, the QoI related to the solution may live in a low-dimensional manifold. We detected this low-dimensionality by exploring the Hessian of the QoI with respect to the parameter and proposed to sample the parameter from a subspace spanned by the eigenvectors of the Hessian instead of the whole parameter space. For the computation of the eigenpairs of the Hessian, we proposed to use a randomized SVD algorithm, whose cost depends only on the number of eigenpairs, not on the nominal dimension of the parameter.
Based on a diffusion model, we demonstrated that this sampling method leads to much smaller errors of the reduced basis approximation of the QoI for parameters with both uniform distribution and Gaussian distribution.

Further investigation includes adaptive Hessian-based sampling on a systematic way to determine the number of eigenvectors and the required accuracy for the reduced basis approximation of the QoI. Another direction is on the development of the Hessian-based sampling method for nonaffine and nonlinear PDE models and quantities of interest. Moreover, study of the comparison of subspace based sampling methods for function approximation and model reduction using local and global Hessian, as well as gradient information is ongoing.



\appendix 

\section{} It is easy to see that $\bE[\bsg^T_{\bar{\bsp}} (\bsp - \bar{\bsp})] = 0$. We only to verify the quadratic term.
Let $(\lambda_k, \bsphi_k)_{k=1}^K$ be the eigenpairs of $\bH_{\bar{\bsp}}$ where $(\bsphi_k)_{k=1}^K$ form a complete orthonormal basis in $\bR^K$ (with completion if the rank of $\bH_{\bar{\bsp}}$ is smaller than $K$). Let $P_K: \bR^K \to \bR^K$ denote a projection operator defined as 
\beq
P_K \bsv = \sum_{k=1}^K  \bsphi_k \bsphi_k^T \bsv, \quad \forall \bsv \in \bR^K,
\eeq
As $(\bsphi_k)_{k=1}^K$ form a complete orthonormal basis in $\bR^K$, we have $P_K \bsv = \bsv$, $\forall \bsv \in \bR^K$. Therefore, by replacing $\bsp-\bar{\bsp} = P_K (\bsp-\bar{\bsp})$, we have
\beq
\begin{split}
\bE\big[(\bsp-\bar{\bsp})^T \bH_{\bar{\bsp}}(\bsp-\bar{\bsp}) \big] &= \sum_{k,l}\bE\big[(\bsp-\bar{\bsp})^T \bsphi_k \bsphi_k^T\bH_{\bar{\bsp}} \bsphi_l \bsphi_l^T(\bsp-\bar{\bsp})\big]\\
& =  \sum_{k}^K \lambda_k \bE\big[(\bsp-\bar{\bsp})^T \bsphi_k \bsphi_k^T(\bsp-\bar{\bsp})\big] \\
& =  \sum_{k=1}^K  \lambda_k \bsphi_k^T \bE\big[(\bsp-\bar{\bsp}) (\bsp-\bar{\bsp})^T \big] \bsphi_k \\
& = \sum_{k = 1}^K \bsphi_k^T \bC \bH_{\bar{\bsp}} \bsphi_k = \text{tr}( \bC \bH_{\bar{\bsp}} ).
\end{split}
\eeq


\section{From local to global Hessian}
\label{app:globalHessian}
Note that the Hessian $\bH_{\bar{\bsp}}$ is local, evaluated at the mean $\bar{\bsp}$, which may fail to characterize the variation of the QoI globally in the parameter space. To deal with this issue,
we propose two schemes for the computation of a global Hessian---namely, an averaged Hessian and a combined Hessian---to account for the variation of the QoI globally in the parameter space.

\subsection{Averaged Hessian}
As the Hessian at the mean describes the local curvature, which might not capture the important directions in the whole parameter spaces, we can replace the Hessian at the mean by an averaged Hessian defined as 
\beq
\bH = \int_P \bH_\bsp d\mu(\bsp) \approx \frac{1}{M} \sum_{m=1}^M \bH_{\bsp^m},
\eeq
with $\bsp^m$ sampled according to its probability distribution $\mu$ in the whole parameter space. Then as in the first case, we can compute the generalized eigenpairs of $(\bH, \bC^{-1})$ with averaged Hessian $\bH$ and project the parameter in the subspace spanned by the first $L$ eigenvectors. We remark that the averaged Hessian was used in \cite{CuiMarzoukWillcox2016} for the projection of a parameter under posterior distribution into a subspace spanned by the dominating eigenvectors in the context of Bayesian inverse problems 

\subsection{Combined Hessian}
\label{sec:combinedHessian}
\begin{figure}[!htb]
\begin{center}
\includegraphics[scale=0.4]{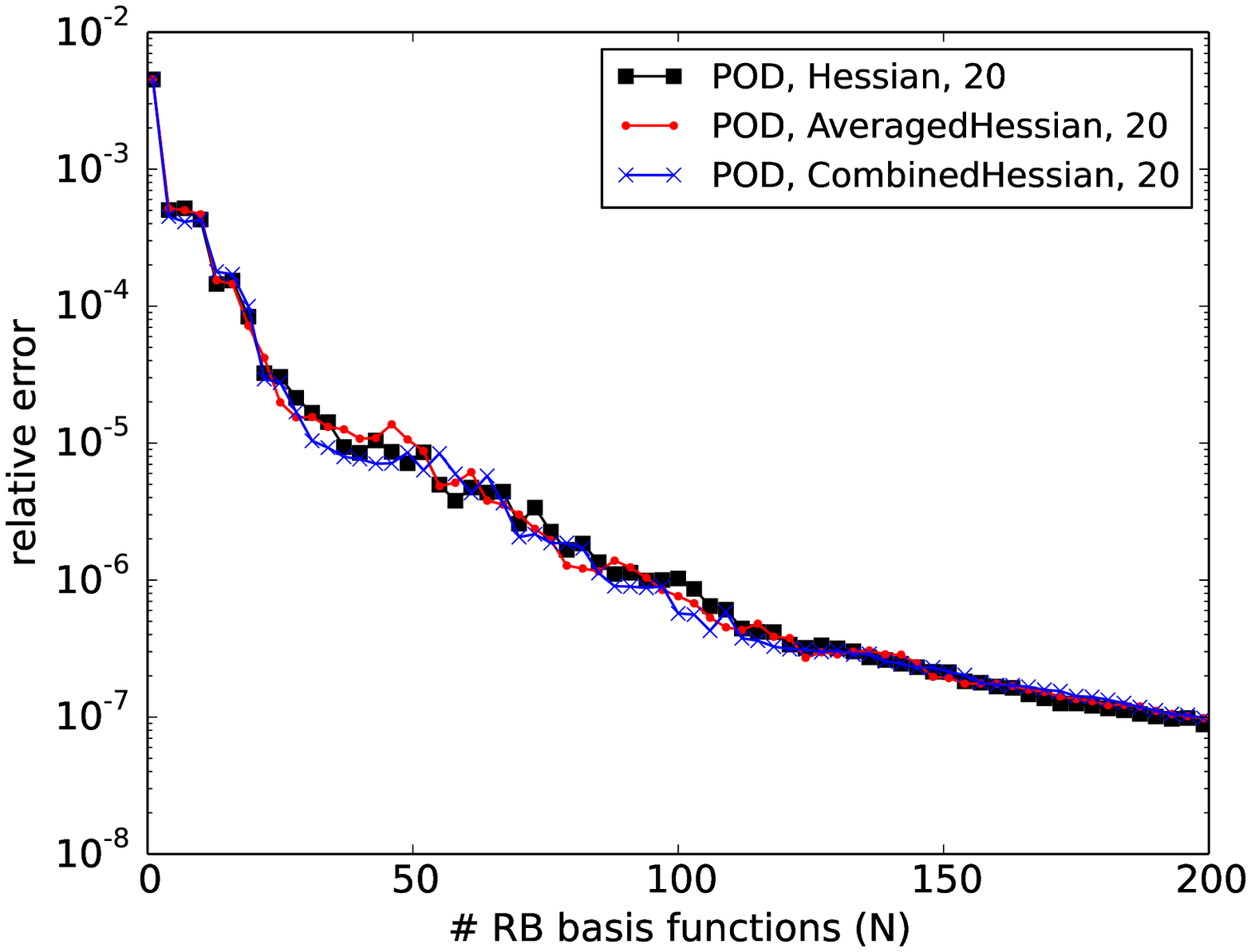}
\includegraphics[scale=0.4]{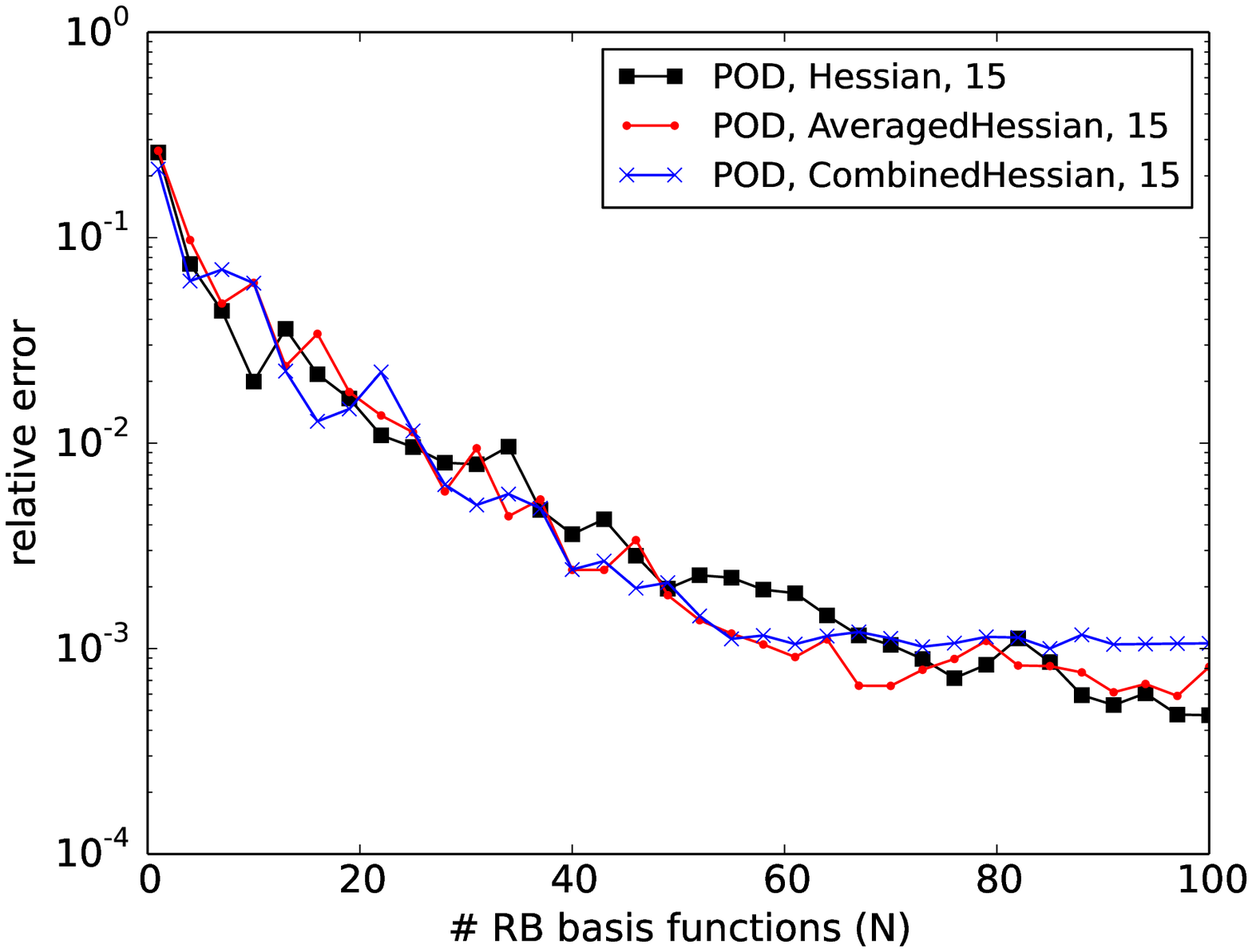}
\end{center}
\caption{Decay of the relative error $\cE_N^s$ defined in \eqref{eq:error} of the RB approximation of the QoI for the Hessian at the mean, the averaged Hessian, and the combined Hessian. Left: example for the uniform distribution in Section \ref{sec:uniform}; right: example for the Gaussian distribution in Section \ref{sec:Gaussian}.
}\label{fig:HessianComparison}
\end{figure}

Another method to construct the subspace for the parameter projection is to combine all the eigenvectors of Hessian at different locations with suitable compression. Let $(\lambda_k^m, \bsvarphi_k^m)_{k=1}^{L_m}$ denote the generalized eigenpairs of $(\bH_{\bsp^m}, \bC^{-1})$ at the random sample $\bsp^m$, $m = 1, \dots, M$, i.e.,
\beq
\bH_{\bsp^m} \bsvarphi_k^m = \lambda_k^m \bC^{-1}\bsvarphi_k^m,  \text{ such that } \bsvarphi_k^m \bC^{-1} \bsvarphi_{k'}^m = \delta_{kk'}, \; k, k' = 1, \dots, L_m, 
\eeq
Then we form the matrix $\Phi \in \bR^{K \times N}$ with $N = L_1 + \cdots + L_M$, by a weighted combination of all the $N$ eigenvectors as 
\beq\label{eq:combinedEigenvector}
\Phi = (w_1^1 \bsvarphi_1^1, \dots, w_{L_1}^1 \bsvarphi_{L_1}^1, \dots, w_1^M \bsvarphi_1^M, \dots, w_{L_M}^M \bsvarphi_{L_M}^M),
\eeq
with suitable weight $w_k^m$ that reflects the importance of the eigenvector $\bsvarphi_k^m$. A natural choice is $ w_k^m =  \sqrt{\lambda_k^m}$. Let $\bS$ denote a matrix such that $\bC^{-1} = \bS \bS^T $, e.g., $\bS = \bC^{-1/2}$ or $\bS$ represents a Cholesky factorization of $\bC^{-1}$. Then we compute the singular value decomposition of $\bS \Phi$ for the largest $L$ singular values as
\beq
\bS\Phi = \bP \Sigma \bQ^T,
\eeq
where $\Sigma = \text{diag}(\sigma_1, \dots, \sigma_L) \in \bR^{L \times L}$ is a diagonal matrix with the $L$ largest non-negative singular values on the diagonal, $\bQ = (\bsq_1, \dots, \bsq_L) \in \bR^{N\times L}$ are the $L$ right singular vectors. We construct the basis functions for the parameter projection as
\beq
\bsvarphi_k = \sum_{j = 1}^{N} \frac{1}{\sigma_k} \bsq_{k,n} \Phi_n, \quad k = 1, \dots, L, 
\eeq
where $\bsq_{k,n}$ denotes the $n$-th element of $\bsq_k$. One can verify that $\bsvarphi_k^T \bC^{-1} \bsvarphi_{k'} = \delta_{kk'}$, $k, k' = 1, \dots, L$.
Fig. \ref{fig:HessianComparison} displays the comparison of the three different Hessians for the construction of the reduced order model, from which we can see that there is almost no difference between using the local Hessian and the global Hessian for the examples in Section \ref{sec:experiments}. Comparison of the differences of the local and global Hessian for both function approximation and model reduction is subject to further investigation.

\section{Multiple quantities of interest}
\label{sec:mulQoI}
In the case of multiple quantities of interest, $s_j, j = 1, \dots, J$, 
instead of constructing different reduced order models with independent Hessian-based sampling for each QoI, we can construct 
a single reduced order model by designing Hessian with the following two approaches as used in Sec. \ref{sec:Hessian-based sampling}. Let $\bH^j_\bsp$ denote the Hessian for $s_j$ at ${\bsp}$, $j =1, \dots, J$, then we can project the full parameter to the eigenvectors of the averaged Hessian at the mean 
\beq
\bH_{\bar{\bsp}} = \frac{1}{J}\sum_{j=1}^J \bH^j_{\bar{\bsp}},  
\eeq
or the double averaged Hessian
\beq
\bH = \frac{1}{JM} \sum_{j=1}^J \sum_{m=1}^M \bH^j_{\bsp^m}.
\eeq
Alternatively, we can compute the eigenpairs of the Hessian $\bH^j_{\bar{\bsp}}$ at mean $\bar{\bsp}$ separately, denoted as $(\lambda_l^j, \bsvarphi_l^j)_{l = 1}^L$,  and combine them with weight $w_l^j$ (e.g. $w_l^j = \sqrt{\lambda_l^j}$ ) as 
\beq
\Phi = (w_1^1 \bsvarphi_1^1, \dots, w_{L_1}^1 \bsvarphi_{L_1}^1, \dots, w_1^J \bsvarphi_1^J, \dots, w_{L_J}^J \bsvarphi_{L_J}^J),
\eeq 
which we compress by SVD as in \ref{sec:combinedHessian} to obtain the dominating singular vectors for the projection of the parameter. Moreover, we may compute the eigenpairs of the sample averaged Hessian with $M$ samples, 
and combine the eigenvectors as
\beq
\Phi = (\Phi_1, \dots, \Phi_J),
\eeq
where $\Phi_j$ is the combined eigenvectors $\eqref{eq:combinedEigenvector}$ for each $j = 1, \dots, J$. 
Then, similarly we perform SVD compression for $\Phi$ and project the full parameter to the subspace formed by the singular vectors corresponding to the dominating singular values.

%
%
%
%
%
%


\end{document}